\documentclass[11pt,reqno]{article}
\usepackage{amssymb,amsmath,amsopn}
\usepackage{amsthm}

\usepackage{setspace}

\newtheorem{theorem}{Theorem}
\newtheorem{lemma}[theorem]{Lemma}

\newtheorem{proposition}[theorem]{Proposition}

\newtheorem*{utheorem}{Theorem}
\newtheorem*{remark}{Remark}

\newcommand{\Cent}{\mathrm{Cent}}

\newcommand{\Sym}{\mathrm{Sym}}
\newcommand{\Alt}{\mathrm{Alt}}
\newcommand{\GF}{\mathrm{GF}}
\newcommand{\GL}{\mathrm{GL}}
\newcommand{\SL}{\mathrm{SL}}
\newcommand{\Mat}{\mathrm{Mat}}

\renewcommand{\deg}{\mathrm{deg\ }}
\renewcommand{\emptyset}{\varnothing}

\newcounter{thmlistcnt}
    {\setcounter{thmlistcnt}{0}%
    \begin{list}{\emph{(\roman{thmlistcnt})}}{%
        \usecounter{thmlistcnt}%
        \setlength{\topsep}{0pt}%
        \setlength{\leftmargin}{0pt}%
        \setlength{\itemsep}{0pt}
        \setlength{\itemindent}{17pt}}%
    }%
    {\end{list}}%

\begin{document}



\setcounter{page}{1}
\title{\bf Commuting conjugacy classes: \\ An application of Hall's Marriage Theorem to Group Theory}
\author{John R. Britnell and Mark Wildon}
\date{25 October 2008}
\maketitle

\thispagestyle{empty}
\begin{center}\bf Abstract\end{center}
Let $G$ be a finite group.
Define a relation $\sim$ on the conjugacy
classes of $G$ by setting $C \sim D$
if there are representatives $c \in C$ and $d \in D$ such that $cd = dc$.
In the case where $G$ has a normal subgroup $H$ such that
$G/H$ is cyclic, two theorems are proved
concerning the distribution, between cosets of $H$, of pairs of conjugacy classes of $G$ related by $\sim$.
One of the proofs involves an interesting application
of the famous Marriage Theorem of Philip Hall.

The paper concludes by discussing some aspects of these theorems and of the 
relation $\sim$
in the particular cases of symmetric and general linear groups,
and by mentioning an open question related to Frobenius groups.

\section{Introduction}
Let~$G$ be a finite group.
Given two conjugacy classes~$C$ and~$D$ of~$G$, we shall say
that~$C$ \emph{commutes} with~$D$, and write~$C \sim D$, if
there exist elements~$c \in C$ and~$d \in D$ such that~$c$ and~$d$
commute.

In this paper we are particularly concerned with the case where~$G$ has
a normal subgroup~$H$ such that~$G/H$ is cyclic; in this case,
since~$G/H$ is abelian, each conjugacy class of~$G$ is entirely
contained within a particular coset of~$H$. We establish some
results concerning the distribution
between the cosets of~$H$ of pairs of conjugacy classes of~$G$
related by~$\sim$.

In order to state our theorems, we make the following definition: a
conjugacy class~$g^G$ of~$G$ is \emph{non-split}
if~$g^G=g^H$, and \emph{split} otherwise.
Observe that~$g^G$ is  non-split if and only if the centre~$\Cent_G(g)$
of~$G$ meets every coset of~$H$ in~$G$, or equivalently, if and only if~$g$
commutes with an element in a generating coset of~$G/H$.
In particular, if a coset~$Ht$ generates the quotient group~$G/H$,
then all of the conjugacy classes in~$Ht$
are non-split.

\begin{theorem}\label{thm:linking}
Let~$G$ be a finite group containing a normal subgroup~$H$
such that~$G/H$ is cyclic.
Let~$Ht$ be a generating coset for~$G/H$,
and let~$Hx$ be any coset.
There is a matching between
the non-split conjugacy classes in~$Hx$
and the conjugacy classes in~$Ht$, such that if~$C$ is matched
with~$D$, then~$C \sim D$.
\end{theorem}
\noindent Our proof of this theorem, given in \S \ref{proofs} below,
involves an interesting application of
Philip Hall's famous Marriage Theorem.

In the special case where~$G/H$
has prime order, Theorem~\ref{thm:linking} can be strengthened
in the following manner.
\begin{theorem}\label{thm:partition}
Let~$G$ be a finite group containing a normal subgroup~$H$
such that~$G/H$ is cyclic of prime order $p$. Let~$Ht$
be a generating coset for~$G/H$.
The non-split conjugacy classes of~$G$ may be partitioned into
sets of the form
\[ \{ g_0^G, g_1^G, \ldots, g_{p-1}^G \} \]
where
$g_m^G \subseteq Ht^m$
and
$g_0^G, g_1^G, \ldots, g_{p-1}^G$ all commute with one another.
\end{theorem}

Theorem~\ref{thm:linking} implies the purely numerical
result that the number of non-split classes in $Hx$
is equal to the number of classes in $Ht$. Indeed, it can be shown that the number of
non-split classes in any two cosets are equal; this is a special case of
\cite[Proposition 9.4]{Isaacs}\footnote{We
thank Tom Wilde for drawing our attention to Isaac's article.}. This fact
has recently been used by the authors \cite{BritnellWildonCosets},
to establish a more general result about the distribution of conjugacy
classes of $G$.

It is not clear whether the assumption in Theorem \ref{thm:linking}, that $Ht$ is a
generating coset, is necessary. The numerical result mentioned in the
preceding paragraph guarantees that for any two cosets~$Hx$ and~$Hy$,
there exists a bijection between the non-split classes of~$Hx$ and the
non-split classes of~$Hy$. We conjecture that a bijection with the
matching property of Theorem~\ref{thm:linking} is always available,
but this remains an open question.

It will be apparent to the reader that we believe
the commuting relation~$\sim$ on conjugacy classes to be of interest
in its own right. In \S \ref{Examples} of
this paper we look at some properties of the relation when~$G$ is a symmetric
group or a finite general linear group. We also observe that if~$|G:H|$ is prime
and all of the non-identity classes in $H$ are split, then~$G$ is a Frobenius group with kernel~$H$.

\section{Proof of Theorems~\ref{thm:linking} and~\ref{thm:partition}}\label{proofs}

We first recall Hall's Marriage Theorem in a form
suitable for our purposes. Hall's original formulation,
and his proof, may be
found in \cite{Hall}.

\begin{utheorem}[Hall's Marriage Theorem]
Suppose that $X$ and $Y$ are finite sets each with~$k$ elements,
and that a relation $\sim$ is defined between~$X$ and~$Y$.
It is possible to order the elements of~$X$ and~$Y$ so that
\begin{align*}
   X &= \left\{ x_1, x_2, \ldots, x_k \right\} \\
   Y &= \left\{ y_1, y_2, \ldots, y_k \right\}
\end{align*}
and
\[ x_i \sim y_i \quad \text{for $1 \le i \le k$} \]
if and only if for every set of~$r$ distinct elements of~$X$, the
total number of elements of~$Y$~which relate to one
of these elements is at least~$r$.
\end{utheorem}

In order to apply Hall's Marriage Theorem to prove Theorem~\ref{thm:linking}
we must show (a)~that the number of non-split conjugacy classes in $Hx$
is equal to the number of conjugacy classes in $Ht$, and (b)~that
given any~$r$ distinct non-split classes from $Hx$, there are at
least~$r$ conjugacy classes in~$Ht$ which commute with
one of these given classes.
The following double
counting argument will establish both of these facts.
We shall use the following notation: if $X$ is any
subset of $G$ and
$g \in G$, then we let $\Cent_X(g) = \Cent_G(g) \cap X$.

Let $C_{i_1}, \ldots, C_{i_r}$ be distinct non-split conjugacy classes in~$Hx$
and let $R = C_{i_1} \cup \ldots \cup C_{i_r}$.
The total
number of conjugacy classes in $Ht$ which commute with one of the
classes making up $R$ is
\[
s =
\sum_{g \in Ht \atop \Cent_R(g) \not= \emptyset}
\frac{1}{|g^G|}
= \frac{1}{|G|}
\sum_{g \in Ht \atop \Cent_R(g) \not= \emptyset}|\Cent_G(g)|. \]
If $g\in Ht$ then $g^G$ is non-split and $\Cent_G(g)$ meets
every coset of $H$ in $G$. In particular, it meets $Hx$.
Using this to rewrite the second sum gives
\[ s =\frac{1}{|H|}
\sum_{g \in Ht \atop \Cent_R(g) \not= \emptyset}|\Cent_{Hx}(g)|.
\]
Now, it is clear that $|\Cent_{Hx}(g)| \ge |\Cent_{R}(g)|$, and so
we can establish a lower bound for $s$ by replacing the condition that $g\in Ht$ with
the condition that $g\in R$ in the sum above. This change makes redundant the second condition, that $\Cent_R(g) \not= \emptyset$; so we have
\[ s \ge \frac{1}{|H|}
\sum_{g \in Ht}|\Cent_{R}(g)|.\]
The quantity on the right-hand side may instead be
found by summing over elements
of~$R$. This gives
\[
\begin{split}
s\ge 
\frac{1}{|H|}\sum_{k \in R} |\Cent_{Ht}(k)| \qquad\qquad\qquad\qquad\qquad\qquad\\\qquad\qquad\qquad\qquad
= \frac{1}{|G|}\sum_{k \in R} |\Cent_{G}(k)| 
= \sum_{k \in R} \frac{1}{|k^G|} = r,
\end{split}
\]
where to pass from the first line
to the second, we use the fact that if $k\in R$ then~$k^G$ is non-split,
and so~$\Cent_{G}(k)$ meets every coset of $H$.
We have shown that~$s\ge r$, which establishes (b). Moreover, if $R$ contains
every non-split class in~$Hx$, then $\Cent_{Hx}(g) = \Cent_R(g)$
for every~$g \in Ht$; hence we have equality at every stage,
and (a) follows. This completes the proof of Theorem~\ref{thm:linking}.\qed
\medskip

\begin{remark}\emph{
It is well known that Hall's Marriage Theorem may be
used to prove that if~$H$ is a subgroup (not necessarily normal) of a
finite group~$G$, then there exists a set of representatives for the
left cosets of $H$ which serves also as a set of representatives for
the right cosets. This fact can be proved without Hall's Marriage
Theorem, or any close equivalent; an elementary proof using double cosets was
given by Miller \cite{Miller} in 1910.}\footnote{Contrary to a folklore
belief, it seems improbable that Hall
was motivated in~\cite{Hall} by the coset
representatives problem, to which we have been unable to find reference
in any of his published work.
It is likely that Hall was aware of Miller's paper and also of Van der
Waerden's paper
\cite{VdV} of 1927. Indeed it is possible that Hall has been confused
at some time with Van der Waerden, who establishes
a theorem similar to Hall's, and explicitly mentions the
problem of coset representatives as his motivation.}
\emph{It seems likely, by contrast, that the use of the Marriage
Theorem in the proof of Theorem~\ref{thm:linking} is essential.}
\end{remark}
\medskip

We now turn to the proof of Theorem~\ref{thm:partition}.
To establish the existence of a partition of the collection of
non-split conjugacy classes of $G$, with the properties stated in the theorem,
we invoke the following lemma.
\begin{lemma}\label{lemma:exponent}
Let~$G$ be a finite group containing a normal subgroup~$H$, and~$Hx$ a coset of $H$.
Let~$c$ be an integer coprime with~$|G|$.
The function on~$G$ given by $g \mapsto g^c$ induces
a bijection between the conjugacy classes in $Hx$ and the
conjugacy classes in~$Hx^c$.
\end{lemma}

\begin{proof}
Since~$c$ is coprime with $|G|$, there exists
an integer~$d$ such that \hbox{$cd \equiv 1$} mod~$|G|$.
It follows that the function $g \mapsto g^c$ is invertible,
with inverse \hbox{$g \mapsto g^d$}.
Moreover, if~$y$ and~$z$ lie in $Hx$ and~$y^c$
and~$z^c$ are conjugate, then~$y^{cd} = y$
and $z^{cd} = z$ are also conjugate. Hence
the function induced on conjugacy classes is bijective.
\end{proof}

Now suppose that in the matching given by Theorem~\ref{thm:linking}
between the non-split conjugacy classes
in~$H$ and the conjugacy classes in~$Ht$,
the class $h^G \subseteq H$
is paired with $g^G \subseteq Ht$
where $h$ and $g$ commute.
For $i = 2, \ldots, p-1$, let $c_i$ be an integer coprime
with~$|G|$ such that $c_i \equiv i \bmod p$.

The conjugacy classes in the set
\[ \bigl\{ h^G, g^G, (g^{c_2})^G, \ldots, (g^{c_{p-1}})^G \bigr\} \]
lie in distinct cosets of $H$
since the integers $c_i$ together with $0$ and $1$ form a complete set of residues modulo $p$. Moreover, the elements~$h$,~$g$ and~$g^{c_i}$ for
$i=2, \ldots, p-1$
certainly commute with one another.
Consider the collection of
sets which are obtained in this way. Lemma~\ref{lemma:exponent} ensures that distinct sets are disjoint,
and it follows easily that they partition the set of non-split classes of $G$. This completes the proof of
Theorem~\ref{thm:partition}.\qed

\section{Examples}\label{Examples}

\subsection{Symmetric groups}
We recall that conjugacy classes of $\Sym(n)$ are parametrized, via their cycle structure, by partitions of $n$.
The classes which lie in $\Alt(n)$ correspond to those partitions for which the number of parts of even size is even.
It is a standard result (see, for example,~\cite[page~65]{FH})
that the class $C^{\lambda}$ corresponding to the partition $\lambda$
splits when the conjugacy action is restricted to $\Alt(n)$,
if and only if the sizes of the parts of $\lambda$ are odd and distinct.

Let $\mathcal{P}(n)$ denote the set of all partitions of $n$. Let $\mathcal{P}_{\textrm{even}}(n)$
be the subset of partitions with an even number of even parts,  $\mathcal{P}_{\textrm{odd}}(n)$ the subset of partitions
with an odd number of even parts, and $\mathcal{D}_{\textrm{o}}(n)$ the subset of partitions
with distinct odd parts. Theorem~\ref{thm:linking} asserts that there exists an bijection
\[ f : \mathcal{P}_{\textrm{even}}(n)\,\backslash\,\mathcal{D}_{\textrm{o}}(n) \longrightarrow
\mathcal{P}_{\textrm{odd}}(n) \]
such that $C^\lambda \sim C^{f(\lambda)}$ for all $\lambda \in \mathcal{P}_{\textrm{even}}(n)\, \backslash\, \mathcal{D}_{\textrm{o}}(n)$.

The numerical implication that $|\mathcal{P}_{\textrm{even}}(n|=|\mathcal{P}_{\textrm{odd}}(n)|-|\mathcal{D}_{\textrm{o}}(n)|$
is well known, but constructing an explicit bijection is
by no means trivial. An elegant one is given by Gupta \cite{Gupta}, which happens to possess the commuting property
in which we are interested.

It is straightforward to describe the commuting relation~$\sim$ in terms of the partitions which parametrize the classes
of $\Sym(n)$. Given two partitions~$\mu$ and ~$\nu$ in $\mathcal{P}(n)$,
we say that~$\nu$ is a \emph{coarsening} of~$\mu$
if~$\nu$ can be obtained from~$\mu$ by
adding together parts of~$\mu$ of the same size.
For example, $(4,3^4,1^2)$ has both $(12,4,2)$
and $(6,4,3^2,1^2)$ as coarsenings.
The following proposition completely describes~$\sim$ for $\Sym(n)$.
\begin{proposition}\label{prop:sym}
Let $\lambda, \mu \in \mathcal{P}(n)$. The conjugacy
classes $C^\lambda$ and $C^\mu$ commute if and only
if there is a partition $\nu \in \mathcal{P}(n)$
which is a coarsening of both $\lambda$ and $\mu$.\hfill$\Box$
\end{proposition}
\begin{proof} The action of a
permutation in~$C^{\nu}$ decomposes into orbits whose lengths are the parts of~$\nu$; a class~$C_{\mu}$ contains an element which acts regularly on
each of these orbits if and only if~$\nu$ is a coarsening of~$\mu$.

Suppose that $x\in C^{\lambda}$ and $y\in C^{\mu}$ are commuting permutations. For $i,j\le n$, let $w_{ij}$ be the subset of $\{1,\dots, n\}$
consisting of elements lying in
$i$-cycles of $x$ and $j$-cycles of $y$. We note that $w_{ij}$ is a union of cycles both of~$x$ and~$y$. If~$\nu$ is the partition whose
parts are the sizes of the sets~$w_{ij}$, then~$\nu$ is a coarsening of both~$\lambda$ and~$\mu$ as required.

Conversely, suppose that~$\nu$
is a coarsening of~$\lambda$ and~$\mu$.
Each part~$k$ of~$\nu$ is an amalgamation of parts of~$\lambda$ with size~$i$, say. If $z\in C^{\nu}$, then~$z^{k/i}$
acts regularly on a $k$-orbit of~$z$, with cycles of length $i$. Thus we may construct an element~$x$ of $C^{\lambda}$
which acts as a power of~$z$ on each orbit of~$z$. Similarly, we can find an element~$y$ of~$C^{\mu}$
which acts as a power of~$z$ on each orbit of~$z$, and clearly~$x$ and~$y$ commute.
\end{proof}

It is worth remarking that Proposition \ref{prop:sym} is essentially about conjugacy classes, rather than individual permutations. In general it is
not the case that if $x,y\in\Sym(n)$ commute then there exists a permutation~$z$, such that~$x$ and~$y$ act as powers of~$z$ on each orbit of~$z$.
The double transpositions in $\Sym(4)$ provide a simple counterexample.

\subsection{Finite general linear groups\protect\footnotemark[3]}
\footnotetext[3]{Certain facts concerning commuting elements in linear groups are stated here
without proof. A paper is in preparation, in which these and other results will be established.}

A natural description of conjugacy classes in general linear groups is provided by the theory of rational canonical form;
this yields the following combinatorial parametrization: for each conjugacy class of $\GL_d(q)$,
a partition $\lambda_f$ is assigned to each monic irreducible polynomial $f$
over $\GF(q)$ other than $f(t)=t$; the only constraint is that
$\sum_{f}|\lambda_f|\cdot\deg f\ =\ d$. If the conjugacy class of an element $M\in\GL_d(q)$ is parametrized with the
assignment~$\{\lambda_f\}$ then each part $a$ of $\lambda_f$ corresponds to an elementary divisor $f^a$ of $M$. The determinant of $M$ is equal to the
product
\[
(-1)^d\prod_f f(0)^{|\lambda|},
\]
from which it is possible to determine, from the parametrization, in which coset of $\SL_d(q)$ a conjugacy class lies.

The commuting relation $\sim$ on conjugacy classes appears to be very much harder to analyse for linear groups than it is for permutation
groups. Working by analogy with Proposition \ref{prop:sym}, it might be tempting to conjecture that if two classes $C_1$ and $C_2$ of $\GL_d(q)$
commute, then there exist elements $Z\in\GL_d(q)$, $X\in C_1$, and $Y\in C_2$, which preserve a common direct sum decomposition, such that $X$ and $Y$
act as polynomials in $Z$ on each summand. However this conjecture would be false; the classes of unipotent elements of $\GL_4(q)$
corresponding to the partitions $(3,1)$ and $(2^2)$ commute, but they have no common proper decomposition,
and no cyclic subalgebra of $\Mat_4(q)$ contains elements of both classes.

A conjugacy class with parametrization $\{\lambda_f\}$ splits when the conjugacy action is restricted to
$\SL_d(q)$, if and only if it satisfies the condition that some prime divisor $d$ of $q-1$ divides all of the parts in all of the
partitions~$\lambda_f$.

Let $\xi$ be a non-zero element of $\GF(q)$ with multiplicative order $q-1$.
Theorem \ref{thm:linking} establishes that there exists a bijection between the non-split classes in $\SL_d(q)$, and the classes
in the coset of $\SL_d(q)$ whose elements have determinant $\xi$. It is difficult to exhibit such a bijection in general, though it is simple
enough to do so in small dimensional cases; we present here the elementary case $d=2$, in order to give an indication of how such bijections can be
constructed. When $q$ is even, the presence of scalar matrices in every coset of
$\SL_2(q)$ means that the problem is trivial, and so we assume that $q$ is odd.

In this case there are four distinct types of conjugacy class: (A) scalars; (B) non-semisimple classes of elements with a single eigenvalue;
(C) classes of elements with two distinct eigenvalues; (D) eigenvalue-free classes. The parameters for (A) and (B) involve only a single linear
polynomial, with the associated partitions being $(1^2)$ and $(2)$ respectively. Classes of type (C) involves two linear polynomials, each with associated
partition~$(1)$. Classes of type (D) involve a quadratic polynomial, with the associated partition~$(1)$.
The split classes, over fields of odd characteristic, are those of type (B).

Let $C_{\xi}$ be the set of elements of $\GL_2(q)$ with determinant $\xi$. It is not hard to calculate the number of classes of each
of the types (A), (B), (C) and (D) in $\SL_d(q)$ and in $C_{\xi}$; the results are presented in the following table.
\vspace{2ex}

\begin{center}
\begin{tabular}{c|cccc}
                    & (A) & (B) & (C)       & (D)       \\ \hline
 $\SL_d(q)$ & $2$ & $2$ & $(q-3)/2$ & $(q-1)/2$ \\
 $C_{\xi}$        & $0$ & $0$ & $(q-1)/2$ & $(q+1)/2$
\end{tabular}
\end{center}
\vspace{2ex}
A useful fact, which we shall not prove here, is that two classes of the same type contain elements which are polynomial in each other, and hence commute.
With this observation, it is easy to present a scheme for a matching between the classes of types (A), (C) and (D) in the two cosets.
We simply match up classes within the types (C) and (D) as far as possible; this leaves just one remaining class of each type in the coset
$C_{\xi}$. Since the two remaining classes in $\SL_d(q)$ contain scalars, the matching can certainly now be completed.

\subsection{Groups of prime index with many split classes}

In the case that $|G:H|$ is prime, there is an interesting characterization of groups with the property that all of the classes in $H$,
other than the identity class, are split.
\begin{proposition}\label{NonIdentityClassesSplit}
Let $G$ be a finite group containing a normal subgroup $H$ of prime index $p$. Every non-identity class of $H$ is split
if and only if $G$ is a Frobenius group with kernel $H$
and complement $C_p$.
\end{proposition}

\begin{proof}
Suppose that every non-identity class of~$H$ is split.
If $t \in G \backslash H$ then~$t$ does not commute with any non-identity element of~$H$.
Since $t^p\in H$, it follows that
$t^p = 1$, and moreover that~$t$ must act on~$H$ as a fixed-point-free
automorphism of order~$p$. It follows that~$G$ is a Frobenius
group with kernel~$H$ and complement~$\left<t\right>$.
The converse follows from standard results on Frobenius groups:
see for example~\cite[Ch.~2, Theorem~7.6]{Gorenstein}.
\end{proof}

It follows from a famous theorem of Thompson (see \cite[Theorem~1]{Thompson} or \cite[Ch.~10, Theorem~2.1]{Gorenstein}) that
the Frobenius kernel $H$ is nilpotent. It is natural to ask what happens when we weaken the hypothesis of Proposition~\ref{NonIdentityClassesSplit}
to require only that every \emph{non-central}
class of $G$ is split. At least when $p=2$, it remains
the case that $H$ must be nilpotent (though not necessarily abelian), but we have no general answer to this question at present.

We conclude with the remark that the centre of $G$ can be recognized by means of the commuting relation $\sim$ on conjugacy classes.
Suppose that $g^G$ is a conjugacy class which commutes with every
other conjugacy class of $G$. Then clearly the conjugates of $\Cent_G(g)$ cover $G$. It is a well-known fact (see for example
\cite[Exercise 237]{Rose})
that the conjugates of a proper subgroup of $G$ cannot cover $G$. Hence $\Cent_G(g)=G$, and so the element $g$ is central.

\end{document}